\documentclass[a4paper,11pt,openright, twosides]{amsart}
\usepackage[utf8]{inputenc}
\usepackage[french]{babel}

\usepackage{amssymb,amsmath,amsfonts}
\usepackage[foot]{amsaddr}

\usepackage{graphicx}
\usepackage{color}
\usepackage{tikz}
\usepackage{natbib}
\usepackage{mathrsfs}           
\usepackage[normalem]{ulem}     

\usepackage{float}
\usepackage{hyperref}

\newcommand{\be}{\begin{equation}}
\newcommand{\ee}{\end{equation}}
\newcommand{\bpm}{\begin{pmatrix}}
\newcommand{\epm}{\end{pmatrix}}

\newcommand{\comment}[1]{}

\newcommand\km{{\mbox{ km}}}
\DeclareTextSymbol{\degre}{OT1}{23}

\begin{document}

\title[La stabilit\'e des lunes de Saturne, Janus et \'Epim\'eth\'ee]{La stabilit\'e des lunes de Saturne,\\Janus et \'Epim\'eth\'ee: de l'observation astronomique \`a la th\'eorie KAM}

\author{Alexandre Pousse$^1$}\address[]{$^1$\textit{Dipartimento di Matematica ``Tullio Levi-Civita"}, Universit\`a degli Studi di Padova, \textit{soutenu par le projet H2020-ERC 677793 StableChaoticPlanetM} }
\author{Laurent Niederman$^{2,3}$}
\address[]{\textit{$^2$Topologie et Dynamique}, LMO, Universit\'e Paris-Saclay, \textit{soutenu par le projet ANR BEKAM (ANR-15-CE40-0001)}}
\author{Philippe Robutel$^3$}
\address[]{$^3$\textit{Astronomie et Syst\`emes Dynamiques}, IMCCE, Observatoire de Paris, Universit\'e PSL}

\date{\today}

\maketitle

\begin{abstract}
Popular science article associated with the work ``On the co-orbital motion in the three-body problem: existence of quasi-periodic horseshoe-shaped orbits" (Arxiv.1806.07262) from the same authors. 
Janus and Epimetheus are two moons of Saturn which exhibit a really peculiar dynamics.
As they orbit on circular trajectories whose radii are only 50 km apart (less than their respective diameters), every four (terrestrial) years the bodies are getting closer and their mutual gravitational influence leads to a swapping of the orbits: the outer moon becoming the inner one and vice-versa.
In this article, we describe how, from this specific astronomical problem to the KAM theory, we came to prove the existence of perpetually stable trajectories associated with the Janus and Epimetheus orbits.

\end{abstract}

\begin{center}
\textsc{Acknowledgment}
\end{center}

Le travail de recherche associ\'e \`a cet article de vulgarisation a \'et\'e soutenu par le projet H2020-ERC 677793 StableChaoticPlanetM.

\tableofcontents

\newpage
\section*{Introduction}

\begin{sloppypar}

La sonde Voyager 1, lanc\'ee en 1977 afin d'explorer Jupiter, Saturne et leurs cort\`eges de satellites, confirma en 1980  l'existence de deux nouvelles lunes de Saturne: Janus et \'Epim\'eth\'ee.
En comparaison avec les lunes Titan, Japet, Rh\'ea, Dione ou T\'ethys dont le diam\`etre d\'epasse le millier de kilom\`etres, Janus et \'Epim\'eth\'ee sont des petits objets avec un diam\`etre de l'ordre de $100\km$.
Toutefois, ces satellites ont la remarquable propri\'et\'e d'\^etre les seuls corps connus du syst\`eme solaire \`a \'evoluer sur des trajectoires en ``fer-\`a-cheval" 
	\footnote{Remarque pour les sp\'ecialistes : ce type d'orbite est sans rapport avec les orbites en fer-\`a-cheval au sens de Smale en dynamique hyperbolique.}
o\`u tous les quatre ans, l'un et l'autre se rapprochent et \'echangent leurs orbites: le plus proche de Saturne passant sur l'orbite externe et inversement (voir la figure \ref{fig:ExchangeOrbits}  et la vid\'eo
	\footnote{\url{https://photojournal.jpl.nasa.gov/catalog/PIA08348}}
 de la NASA issue des observations de la sonde Cassini).

Il est naturel de se demander si cette configuration observ\'ee depuis plus de 30 ans est stable sur des temps beaucoup plus longs.
Autrement dit: est-ce que Janus et \'Epim\'eth\'ee auront des rapprochements et cons\'ecutivement des \'echanges d'orbites sur des temps comparables \`a l'\^age du syst\`eme solaire?
Ou bien, est-ce qu'au bout d'un certain temps ces \'echanges n'auront plus lieu \`a cause, par exemple, d'une collision mutuelle ou, au contraire d'un \'eloignement progressif des deux lunes?

Dans un travail r\'ecent, nous avons \'etudi\'e une id\'ealisation de ce syst\`eme (Saturne, Janus et \'Epim\'eth\'ee repr\'esent\'ees par trois corps massifs ponctuels \'evoluant dans un m\^eme plan) et d\'emontr\'e rigoureusement (en un sens qui sera pr\'ecis\'e par la suite) l'existence de mouvements  perp\'etuellement stables de ce type.
Plus pr\'ecis\'ement, \textbf{\textit{notre r\'esultat garantit l'existence de conditions initiales (positions et vitesses des deux lunes \`a un instant donn\'e) telles que le ballet gravitationnel men\'e par Janus et \'Epim\'eth\'ee continuera ind\'efiniment.}}

Le point de d\'epart de ce travail est un mod\`ele d\'evelopp\'e par P. Robutel et A. Pousse pour l'\'etude de la dynamique de deux plan\`etes en r\'esonance co-orbitale (i.e. dont les p\'eriodes orbitales sont identiques; plus bas, dans l'encadr\'e du premier paragraphe nous revenons sur ce ph\'enom\`ene de r\'esonance).
Les motivations \'etaient purement astronomiques, li\'ees aux recherches d'exoplan\`etes (donc de plan\`etes en dehors de notre syst\`eme solaire) dans cette configuration particuli\`ere.  
En effet, parmi les milliers de d\'etections accumul\'ees depuis 1995 et malgr\'e une grande diversit\'e de dynamiques observ\'ees, pour l'instant aucun syst\`eme de plan\`etes extrasolaires ne comporte de couple de plan\`etes ``co-orbitantes".
Or, a priori aucun ph\'enom\`ene physique n'entrave leur formation et, au contraire, certaines \'etudes num\'eriques (c'est-\`a-dire des calculs de solutions r\'ealis\'es par ordinateur) sugg\`erent leur existence.
D'un point de vue plus th\'eorique, il s'agissait donc de construire un mod\`ele  permettant d'interpr\'eter les exp\'erimentations num\'eriques qui avaient \'et\'e men\'ees pour ce type de probl\`eme.
Puis, lors d'une discussion de couloir, nous avons r\'ealis\'e qu'une d\'emonstration rigoureuse de la stabilit\'e perp\'etuelle de mouvements co-orbitaux de m\^eme type que ceux de Janus et \'Epim\'eth\'ee  pouvait \^etre obtenue en associant \`a ce mod\`ele des th\'eor\`emes de stabilit\'e d\'evelopp\'es par J. P\"oschel
\footnote{Remarque: ce sont des th\'eor\`emes  qui ont \'et\'e obtenus initialement par S. Kuksin et H. Eliasson.}
 dans les ann\'ees 90.

\medskip

Apr\`es une pr\'esentation d\'etaill\'ee de la probl\'ematique astronomique pos\'ee par les trajectoires de Janus et \'Epim\'eth\'ee (voir le premier paragraphe), nous aborderons le th\`eme plus g\'en\'eral de la stabilit\'e en m\'ecanique c\'eleste (second paragraphe).
C'est une question tr\`es ancienne qui remonte \`a Newton avec la d\'ecouverte de l'attraction universelle.
Il est g\'en\'eralement tr\`es difficile d'y r\'epondre m\^eme lorsque le syst\`eme physique consid\'er\'e est restreint \`a des corps massifs ponctuels, c'est-\`a-dire au probl\`eme id\'eal et math\'ematiquement bien d\'efini connu sous le nom de ``probl\`eme des $N$ corps". 

Alors que le cas $N=2$ a \'et\'e enti\`erement r\'esolu par Newton,
pour $N=3$, seuls quelques rares r\'esultats de stabilit\'e perp\'etuelle ont \'et\'e d\'emontr\'es, notamment dans le cas des configurations d'\'equilibre d'Euler et de Lagrange (le troisi\`eme paragraphe est d\'edi\'e \`a ce sujet) et pour la configuration Soleil-Jupiter-Saturne dans un article
	\footnote{Plus pr\'ecis\'ement, le r\'esultat \'etait valable pour trois corps \'evoluant dans un m\^eme plan et a \'et\'e \'etendu dans le cas g\'en\'eral ($N$ corps dans l'espace) par M. Herman et J. F\'ejoz.} 
fondateur d'Arnol'd en 1963 (voir le cinqui\`eme paragraphe).

La situation \`a consid\'erer pour les orbites en fer-\`a-cheval est particuli\`erement complexe (les d\'etails sont pr\'esent\'es dans le quatri\`eme paragraphe) ce qui a rendu difficile l'application des th\'eories existantes bien que l'existence de ces trajectoires ait \'et\'e pressentie au d\'ebut du vingti\`eme si\`ecle par l'astronome E. Brown 
	\footnote{Intuition remarquable puisque \`a l'\'epoque aucun corps c\'eleste en dynamique ``fer-\`a-cheval" n'avait \'et\'e observ\'e.}.
C'est en poursuivant l'id\'ee de l'article d'Arnol'd tout en surmontant les obstacles li\'es \`a notre situation compliqu\'ee  \`a l'aide de th\'eor\`emes beaucoup plus r\'ecents (voir le dernier paragraphe) que nous avons pu \'etablir notre r\'esultat de stabilit\'e perp\'etuelle pour des trajectoires de type fer-\`a-cheval dans le probl\`eme des 3 corps.

\medskip

D'un point de vue purement astronomique, notre travail n'apporte pas de r\'eponse d\'efinitive quand \`a la stabilit\'e d'\'eventuelles exoplan\`etes en dynamique ``fer-\`a-cheval" ou des trajectoires de Janus et \'Epim\'eth\'ee.
Pour ces derni\`eres, il faudrait tenir compte de nombreux effets suppl\'ementaires comme l'influence gravitationnelle d'autres satellites de Saturne, l'effet des anneaux, l'aplatissement de Saturne, les effets de mar\'ees, etc...
Cependant, notre r\'esultat indique que les r\'eguli\`eres rencontres proches entre les deux lunes ne d\'estabiliseront pas n\'ecessairement leur mouvement.
Pour aller plus loin, il serait n\'ecessaire de raffiner la repr\'esentation du ph\'enom\`ene en int\'egrant les divers effets physiques.
La probl\'ematique deviendrait alors purement astronomique.

Nous faisons remarquer que ce type de travail rapprochant math\'ematiciens et astronomes  n'est pas si courant.
D'une part, c'est un mod\`ele d\'evelopp\'e pour des raisons astronomiques et permettant aussi d'obtenir des r\'esultats math\'ematiques qui a favoris\'e notre association.
D'autre part, il n'est pas simple d'avoir une collaboration interdisciplinaire aussi avanc\'ee, ne serait-ce que par manque de langage commun.

\end{sloppypar}

\section{Janus et \'Epim\'eth\'ee: \'echanges d'orbites et ``fers-\`a-cheval"}

\begin{figure}[H]
	\begin{center}
		\includegraphics[scale=0.575]{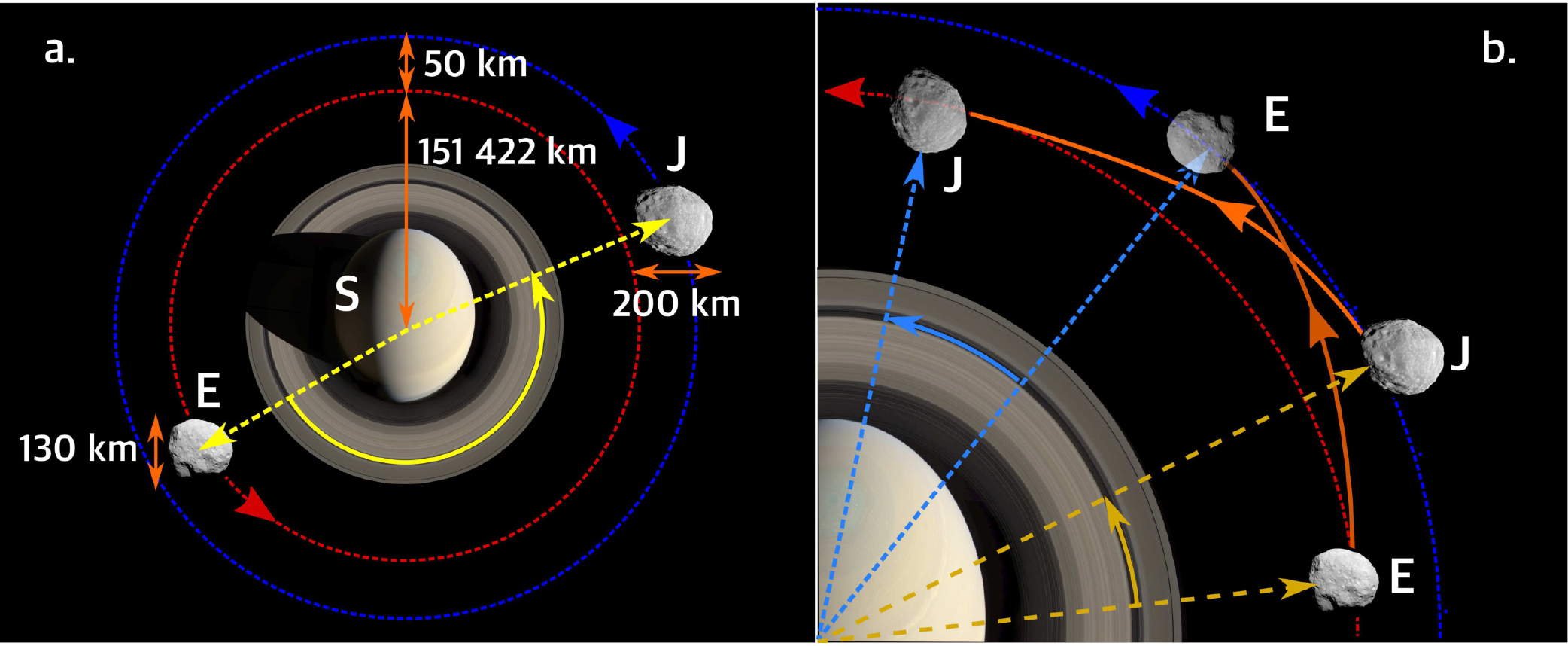}\\[0.25cm]
		\caption{\small{
			(a.) Syst\`eme Saturne-Janus-\'Epim\'eth\'ee lors de l'observation de Voyager 1 en f\'evrier 1980. Les trois corps sont presque align\'es, les lunes de part et d'autre de Saturne, tels qu'\'Epim\'eth\'ee se situe sur l'orbite interne (orbite rouge) tandis que  Janus est sur l'orbite externe (orbite bleue).
			(b.) Repr\'esentation de l'\'echange d'orbite de 1982. Au moment de leur rencontre proche, sans jamais \^etre d\'epass\'ee par \'Epim\'eth\'ee, Janus acc\'el\`ere et se d\'ecale vers l'orbite interne (orbite rouge), tandis qu'\'Epim\'eth\'ee ralentit et ``tombe" sur l'orbite externe (orbite bleue). (\textit{Images: A.Pousse}).
		}}
		\label{fig:SchemaEO}
	\end{center}
\end{figure}

\begin{sloppypar}

Janus ($\mbox{\bf{J}}$) et \'Epim\'eth\'ee ($\mbox{\bf{E}}$) orbitent \`a une distance moyenne de Saturne ($\mbox{\bf{S}}$) de $151~440\km$ (soit $2.5$ fois le rayon de la plan\`ete) en effectuant leur r\'evolution en un peu moins de $17$ heures.
En premi\`ere approximation, ces deux satellites d\'ecrivent dans un m\^eme plan des orbites circulaires dont les rayons diff\`erent d'une cinquantaine de kilom\`etres, c'est-\`a-dire moins que leurs tailles respectives.  
Des collisions sont donc possibles a priori.

En f\'evrier 1980, lors de l'observation de Voyager 1, les deux lunes se situaient de part et d'autre de Saturne, Janus gravitant sur l'orbite la plus externe (voir la figure \ref{fig:SchemaEO}.a).
Une analyse grossi\`ere bas\'ee uniquement sur la 
troisi\`eme
	\footnote{\url{https://fr.wikipedia.org/wiki/Lois_de_Kepler}}
loi de Kepler implique qu'\'Epim\'eth\'ee rattrape lentement Janus ce qui, compte tenu de leurs dimensions, engendrerait in\'evitablement une collision au cours de l'ann\'ee 1982.
Or cette collision n'eut jamais lieu.
En effet, lorsqu'\'Epim\'eth\'ee rattrapa Janus, il y eut \textit{rencontre proche}: leurs interactions gravitationnelles mutuelles devinrent suffisamment fortes pour modifier significativement leurs trajectoires circulaires.
De ce fait, sans jamais \^etre d\'epass\'ee par \'Epim\'eth\'ee, Janus acc\'el\'era en se d\'ecalant vers l'orbite interne tandis qu'\'Epim\'eth\'ee ralent\^it et ``tomba" sur l'orbite externe (voir la figure \ref{fig:SchemaEO}.b).
Janus s'\'eloignant d'\'Epim\'eth\'ee, leurs interactions gravitationnelles mutuelles faiblirent ce qui ``figea" leurs orbites sur des cercles (quasi) identiques \`a l'instant pr\'ec\'edent la rencontre proche.
Il y eut ainsi \textit{\'echange d'orbites} sans que les deux lunes n'entrent en collision, leur distance minimale ayant \'et\'e de $14~000\km$ (l'influence gravitationnelle de Saturne sur la dynamique des lunes restant dominante par rapport aux interactions mutuelles des deux lunes).
Cet  \'echange ne fut pas instantan\'e mais s'effectua de mani\`ere continue et relativement lente pendant environ 8 mois.

\end{sloppypar}

\begin{figure}[H]
	\begin{center}
		\includegraphics[scale=0.1725]{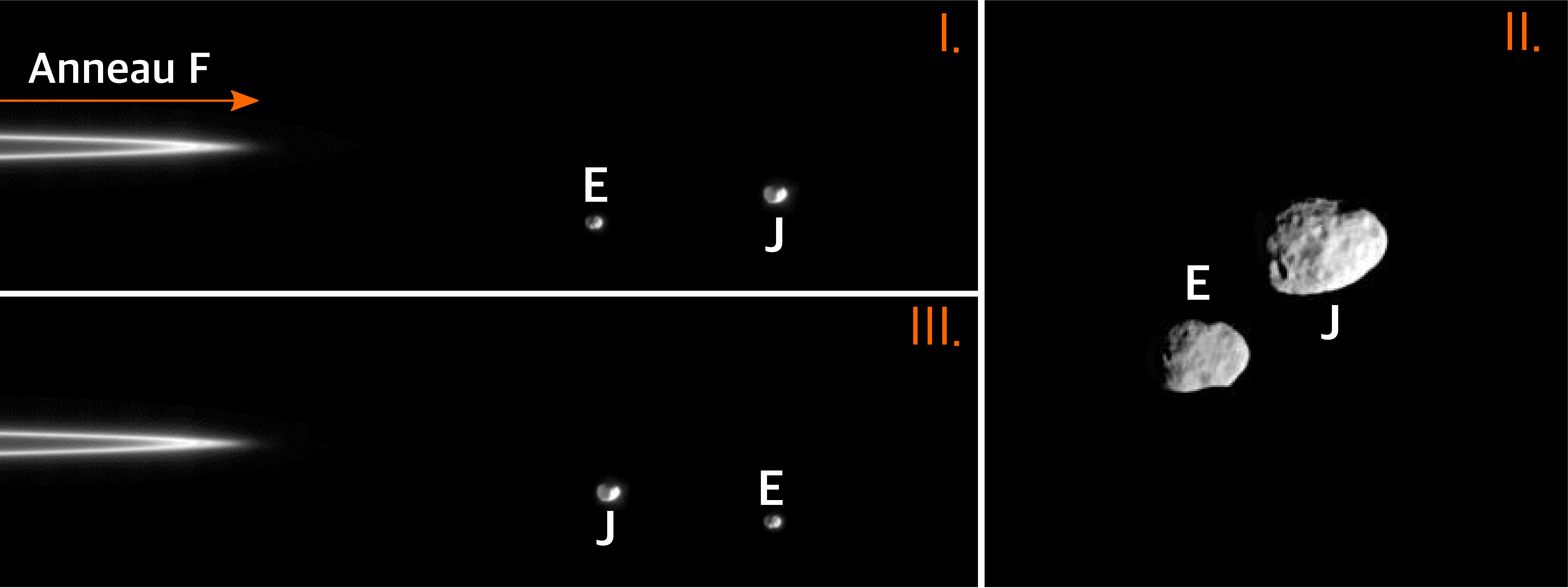}\\[0.5cm]
		\caption{\small{
I-III. Observations par la sonde Cassini en 2006 d'une rencontre proche suivit d'un \'echange d'orbites entre Janus et \'Epim\'eth\'ee. \textit{(Images: NASA/JPL/Space Science Institute)}
		}}
	\label{fig:ExchangeOrbits}
	\end{center}
\end{figure}

\begin{sloppypar}

Depuis sa d\'ecouverte, ce ``ballet" gravitationnel est observ\'e tous les quatre ans (aussi bien par des sondes spatiales orbitant autour de Saturne que par des t\'elescopes professionnels sur Terre). 
La figure \ref{fig:ExchangeOrbits} constitu\'ee de trois images prises par la sonde Cassini en 2006, ainsi qu'une vid\'eo
	\footnote{\url{https://photojournal.jpl.nasa.gov/catalog/PIA08348}} 
de la NASA fournissent l'exemple d'un moment o\`u Janus ``passe devant" \'Epim\'eth\'ee et devient alors la lune interne du syst\`eme (comme observ\'e en 1982; voir la figure \ref{fig:SchemaEO}.b).
Bien que les deux lunes semblent tr\`es proches dans la seconde image, Janus se situe \`a plus de $40~000\km$ d'\'Epim\'eth\'ee.

Le  fait que les rayons des orbites soient presqu'\'egaux et donc que leur p\'eriodes orbitales soient quasi-identiques engendre des comportements particuliers, caract\'eristiques d'une dynamique en \textit{r\'esonance} (comme celle qui appara\^it dans le mouvement entretenu d'une balan\c{c}oire).

\medskip

\noindent\fbox{\parbox{\linewidth\fboxrule\fboxsep}{\textbf{Sur la notion de r\'esonance.}\\
Plus pr\'ecis\'ement, supposons que les orbites circulaires que deux lunes suivraient autour d'une plan\`ete si leurs masses \'etaient infiniment petites soient telles que leurs p\'eriodes orbitales $T$ et $T'$ sont commensurables, c'est-\`a-dire dans un rapport rationnel donc $T/T' = p/q$ o\`u $p,q$ sont des entiers quelconques.
Cela signifieraient que si la plan\`ete et ses deux satellites sont initialement en conjonction alors tous les $q$ tours de la lune $1$ ou les $p$ tours de la lune $2$, les trois corps seraient \`a nouveau dans cette m\^eme configuration de conjonction.
Dans cette situation, les deux lunes sont dites en r\'esonance.

Consid\'erons maintenant une situation \`a proximit\'e de la r\'esonance.
La perturbation provoqu\'ee par les interactions gravitationnelles mutuelles des deux lunes est alors presque p\'eriodique ce qui peut induire une augmentation importante de l'amplitude des orbites et conduire \`a une \'ejection ou une collision.
Toutefois, \`a l'int\'erieur de la ``zone" de r\'esonance, certaines trajectoires peuvent \^etre stables avec des interactions qui se produisent pendant un temps caract\'eristique beaucoup plus long que ce qui l'engendre (un exemple de ce ph\'enom\`ene est mis en \'evidence dans la vid\'eo 
\url{http://images.math.cnrs.fr/Du-ressort-a-l-atome-une-histoire-de-resonance.html}
d'Images des Math\'ematiques r\'ealis\'ee par Beno\^it Grebert et qui est d\'edi\'ee aux r\'esonances dans le cas d'un syst\`eme physique beaucoup plus simple: une corde et deux pinces-\`a-linge).
Dans le cas de Janus et \'Epim\'eth\'ee, le ph\'enom\`ene d'\'echange d'orbites intervient tous les quatre ans alors que leur p\'eriodes orbitales sont de $17$ heures.}}

\medskip

\begin{figure}[H]
	\begin{center}
		\includegraphics[scale=0.575]{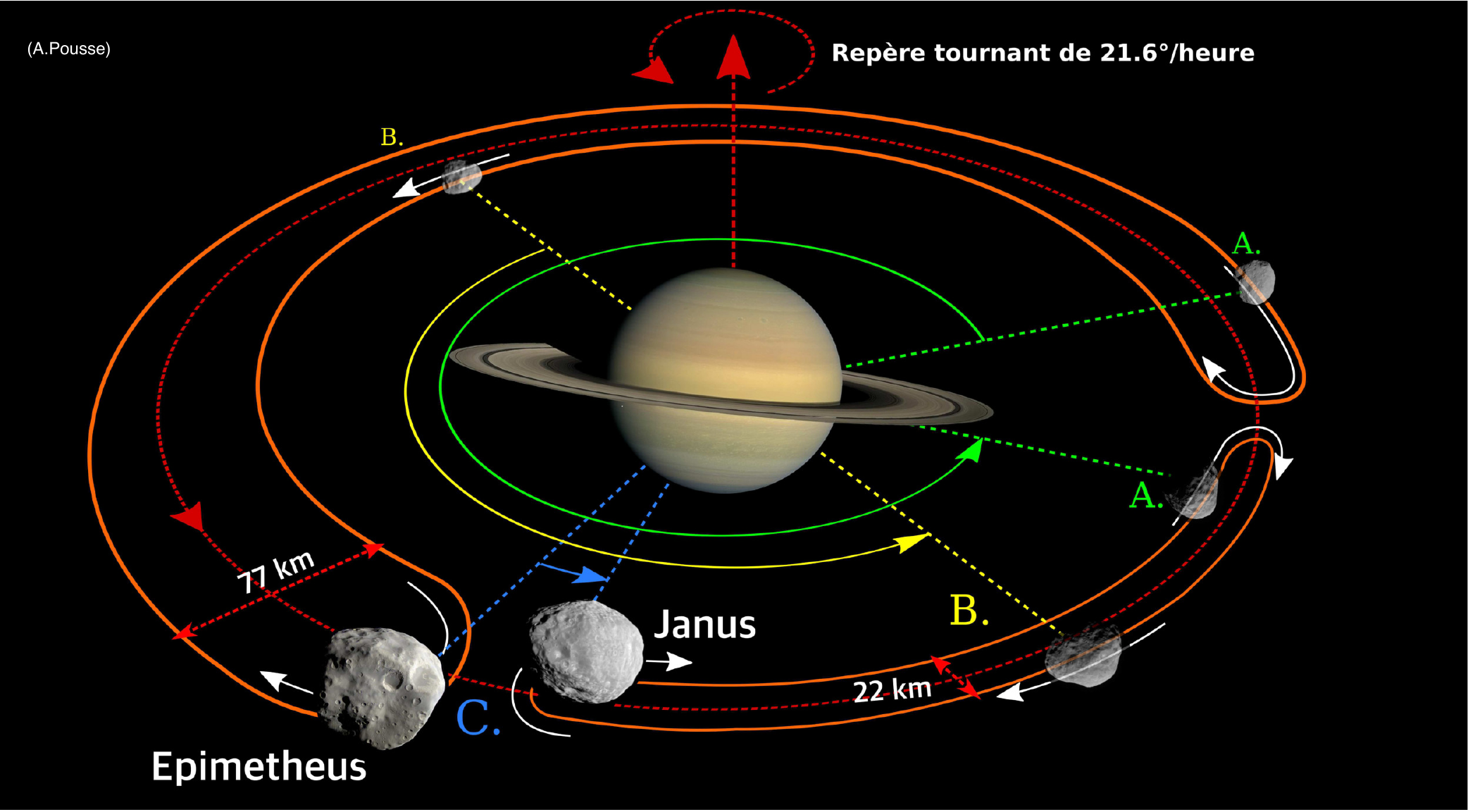}\\[0.5cm]
		\caption{\small{
			Repr\'esentation des trajectoires en fer-\`a-cheval parcourues par Janus et \'Epim\'eth\'ee sur une p\'eriode de 8 ans.
			La configuration \textbf{A} correspond \`a un \'echange d'orbites o\`u \'Epim\'eth\'ee et Janus deviennent respectivement les lunes interne et externe. Les configurations \textbf{B} (trois corps align\'ees tels que les lunes sont de part et d'autre de Saturne) et \textbf{C} (inverse de la configuration \textbf{A}, rencontre proche o\`u \'Epim\'eth\'ee et Janus deviennent respectivement les lunes externe et interne) sont associ\'ees respectivement aux positions des deux lunes lors de l'observation de Voyager 1 en 1980 (figure \ref{fig:SchemaEO}.a) et des \'echanges d'orbites de 1982 (figure \ref{fig:SchemaEO}.b) et de 2006 (figure \ref{fig:ExchangeOrbits}). (\textit{Image: A.Pousse}).
			}}

		\label{fig:SJE}
	\end{center}
\end{figure}

Pour Janus et \'Epim\'eth\'ee, c'est une r\'esonance ``$1/1$" ou ``co-orbitale" qui est \`a l'oeuvre.
Une bonne mani\`ere  d'observer leurs interactions est de repr\'esenter le mouvement relatif des lunes vu par un observateur Saturnien tournant sur lui-m\^eme \`a la vitesse angulaire moyenne de Janus et \'Epim\'eth\'ee.
Ainsi, avec ce r\'ef\'erentiel repr\'esent\'e dans la figure \ref{fig:SJE}, les lunes parcourent des trajectoires en \textit{fer-\`a-cheval} sur une p\'eriode de $8$ ans ce qui correspond \`a plus de $4~000$ r\'evolutions autour de Saturne.

\end{sloppypar}

\section{Newton et la stabilit\'e du syst\`eme solaire}

\begin{sloppypar}
Kepler fut le premier \`a donner une description pr\'ecise du mouvement des plan\`etes suivant ses trois c\'el\`ebres lois
	\footnote{\url{https://fr.wikipedia.org/wiki/Lois_de_Kepler}}. 
Cependant, en enserrant leur mouvement dans un embo\^itement de solides platoniciens
	\footnote{J. Kepler \textit{``Mysterium Cosmographicum"}, 1596; voir \url{https://en.wikipedia.org/wiki/Mysterium_Cosmographicum}},
 il consid\'erait un ``Monde" rigide et immuable o\`u la question de la stabilit\'e ne se posait pas.

Newton r\'evolutionna cette description avec une repr\'esentation dynamique du syst\`eme solaire, r\'egie par un principe fondamental ``\textit{la masse du corps multipli\'e par son acc\'el\'eration est \'egale \`a la somme des forces qui s'exercent sur lui}" associ\'e \`a sa loi d'attraction universelle ``\textit{deux corps s'attirent en raison directe de leur masse et en raison inverse du carr\'e de leur distance}".
L'invention du calcul diff\'erentiel lui permit de r\'eduire le mouvement des corps qui composent un syst\`eme physique \`a un jeu d'\'equations diff\'erentielles
	\footnote{\url{https://fr.wikipedia.org/wiki/Equation_differentielle}} 
du second ordre n\'ecessitant la connaissance de deux quantit\'es, les vitesses et positions \`a un instant initial (ce que l'on d\'enomme plus simplement \textit{conditions initiales}), pour d\'eterminer une \textit{solution}: les vitesses et positions au cours du temps.

L'application des lois de Newton au mouvement d'une unique plan\`ete  autour du Soleil bouleversa notre compr\'ehension du monde. 
En effet,  ce \textit{probl\`eme \`a 2 corps} poss\'edant -- miraculeusement -- une solution explicite pour toutes conditions initiales, il permit de retrouver les trois lois \'enonc\'ees par Kepler pour d\'ecrire les trajectoires des plan\`etes du syst\`eme solaire.
Or, si les plan\`etes sont attir\'ees par le Soleil, la loi d'attraction universelle impose aussi qu'elles s'influencent mutuellement.
De ce fait, rien ne justifie qu'elles d\'ecrivent perp\'etuellement un \textit{mouvement keplerien}
	\footnote{\url{https://fr.wikipedia.org/wiki/Mouvement_keplerien}}
comme dans le cas de deux corps : au bout d'un certain temps, collisions, \'ejections, r\'e-agencements deviennent possibles.
Pour Newton, la stabilit\'e du monde telle que nous l'observons ne doit alors son salut qu'\`a l'intervention de temps \`a autre d'un ``Grand horloger"
	\footnote{I. Newton ``\textit{Trait\'e d'optique}", 1732, p.489}. 
Ses d\'ecouvertes firent cependant \'emerger la question de la stabilit\'e du syst\`eme solaire et plus g\'en\'eralement la recherche de solutions dans le \textit{probl\`eme des N corps},
 o\`u $N$ est un entier quelconque, ce qui donna naissance \`a la \textit{m\'ecanique c\'eleste}.

\end{sloppypar}
\begin{figure}[H]
	\begin{center}
		\includegraphics[scale=0.575]{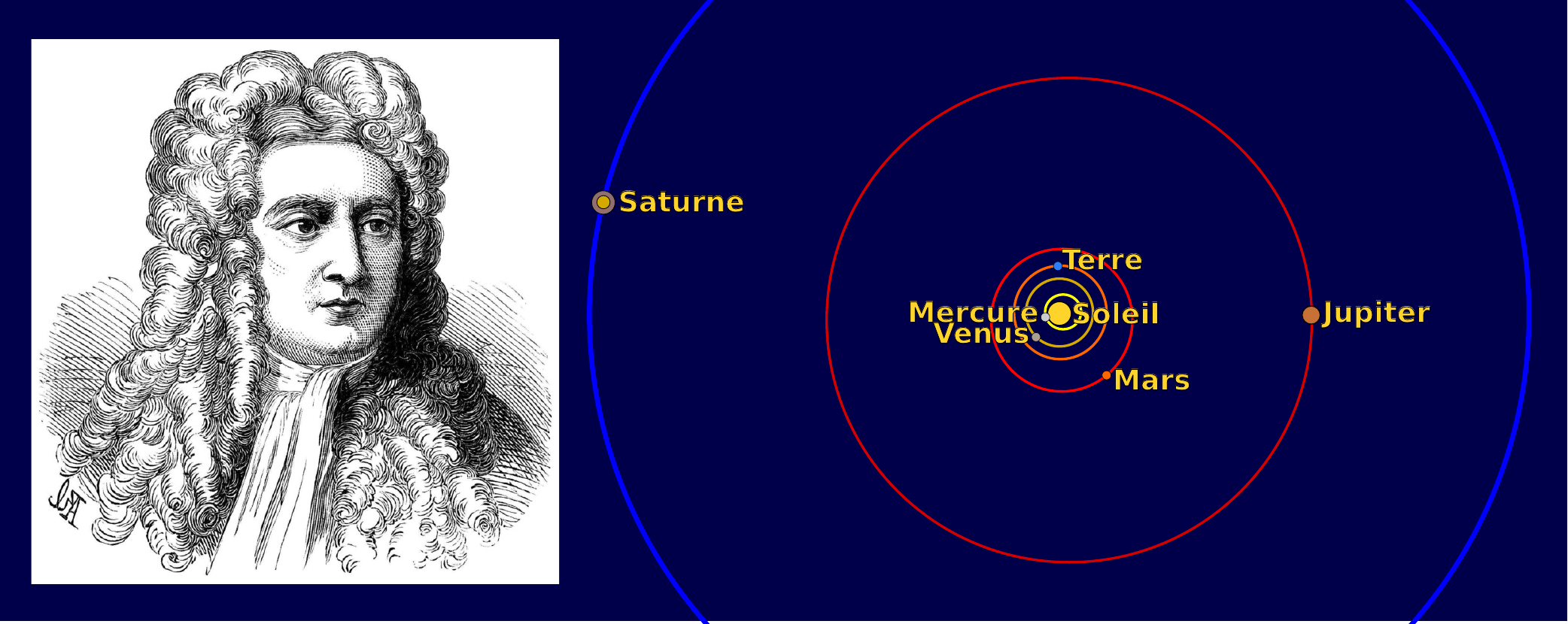}\\[0.25cm]
	\caption{\small{
		Newton et le syst\`eme solaire.
	}}
	\end{center}
\end{figure}

\begin{sloppypar}

Clairaut fut l'un des premiers \`a \'ecrire pr\'ecis\'ement les \'equations du \textit{probl\`eme des 3 corps} et comprendre l'\'etendue des difficult\'es pos\'ees par la recherche de solution. 
D'ailleurs, il termina son calcul par:\textit{``Int\`egre maintenant qui pourra"} 
	\footnote{A. Clairaut ``\textit{M\'emoire lu \`a l'Acad\'emie des Sciences le 23 juin 1759. Contenant des r\'eflexions sur le Probl\`eme des trois corps, avec les \'equations diff\'erentielles qui expriment les conditions de ce Probl\`eme}", 1759, Journal des S\c{c}avans, p.563.}.
Face \`a ce constat, les travaux se sont alors focalis\'es sur la recherche de \textit{solutions d'\'equilibre} (i.e. configurations o\`u le mouvement relatif des corps est fig\'e au cours du temps, voir un exemple plus bas) afin de tenter d'\'eclaircir le probl\`eme des 3 corps.\\

\noindent\fbox{\parbox{\linewidth\fboxrule\fboxsep}{\textbf{Pour approfondir...}\\
Pour plus de d\'etails, nous renvoyons le lecteur \`a la belle exposition virtuelle ``\textit{les math\'ematiques du ciel}" (\url{http://ciel.mmi-lyon.fr}) issue d'une collaboration entre le Labo junior de l'ENS Lyon et le mus\'ee des confluences.
Plus particuli\`erement, la partie ``\textit{2 astres en t\^ete \`a t\^ete}" (\url{http://ciel.mmi-lyon.fr/deux-astres-en-tete-a-tete/}) pr\'esente le probl\`eme des 2 corps (et pourquoi l'\'etudier), les lois et concepts de Kepler, la r\'evolution apport\'ee par Newton et enfin un entretien \'eclairant avec Alain Chenciner \`a propos du probl\`eme des 3 corps.}}

\medskip
\end{sloppypar}

\section{La recherche de solutions dans le probl\`eme des 3 corps}

\begin{figure}[H]
	\begin{center}
		\includegraphics[scale=0.575]{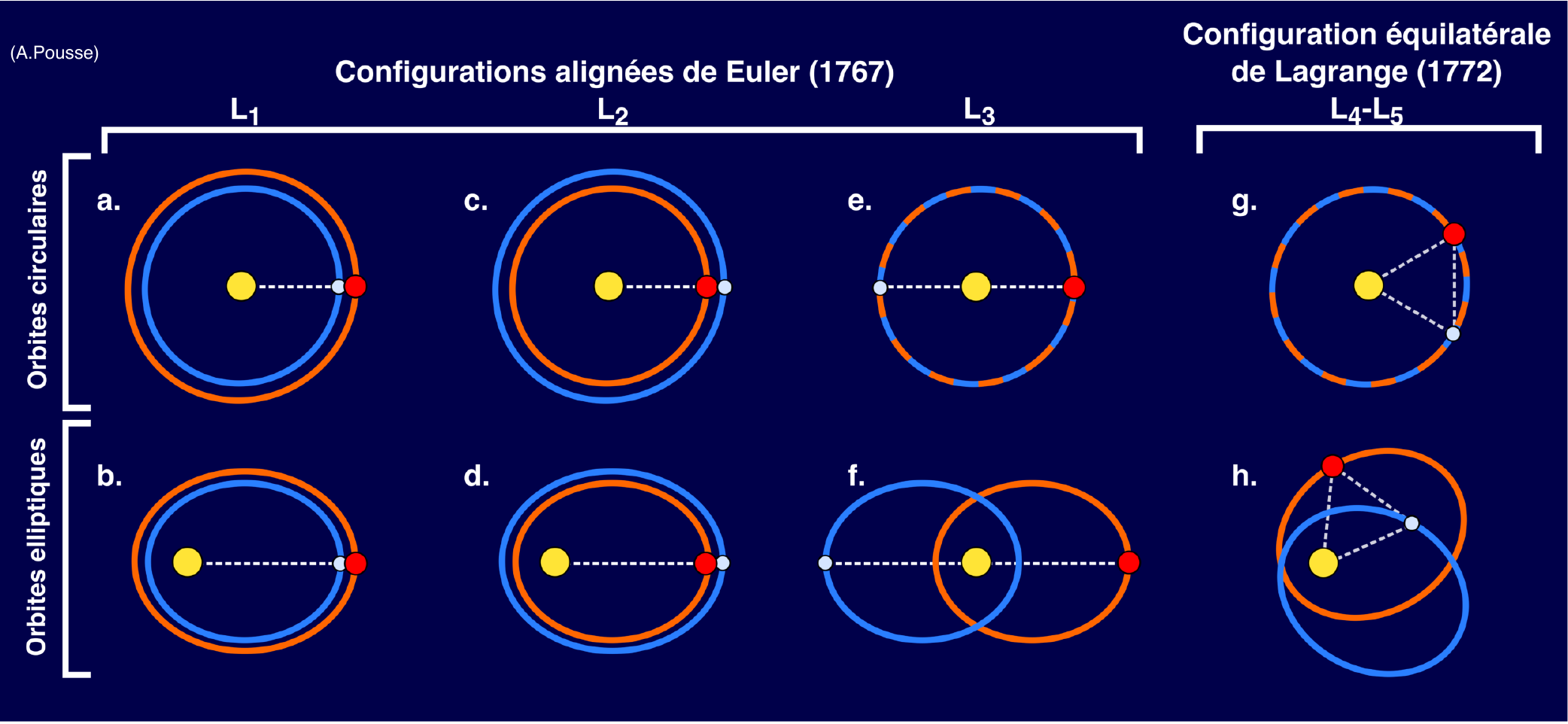}\\[0.25cm]
		\caption{\small{
			Configurations d'\'equilibre dans le probl\`eme des 3 corps: les trois familles de configurations align\'ees d'Euler (a--f) et  la famille de configurations \'equilat\'erales de Lagrange (g--h).}(\textit{Image: A.Pousse}).
		}
	\label{fig:EuLa}
	\end{center}
\end{figure}

\begin{sloppypar}
La premi\`ere solution remarquable du probl\`eme des 3 corps fut pr\'esent\'ee par Euler en 1764 quand il fit remarquer que, \textit{``si la Lune  \'etait quatre fois plus \'eloign\'ee de la Terre qu'elle ne l'est actuellement, les mouvements relatifs du Soleil, de la Terre et de la Lune seraient tels que cette derni\`ere nous appara\^itrait comme une \'eternelle pleine Lune"}.
Plus g\'en\'eralement, il mit en \'evidence trois familles de solutions pour lesquelles deux corps gravitent autour d'un troisi\`eme tels que les trois corps restent en permanence align\'es sur une droite mobile (voir les figures \ref{fig:EuLa}.a--f).
La seconde solution remarquable fut d\'ecouverte par Lagrange en 1772: la configuration \'equilat\'erale. 
Comme cela est illustr\'e sur les figures \ref{fig:EuLa}.g--h, cette famille de solutions repr\'esente une configuration o\`u deux corps gravitent autour d'un troisi\`eme telle que les trois corps forment en permanence un triangle \'equilat\'eral en rotation.

\begin{figure}[H]
	\begin{center}
		\includegraphics[scale=0.575]{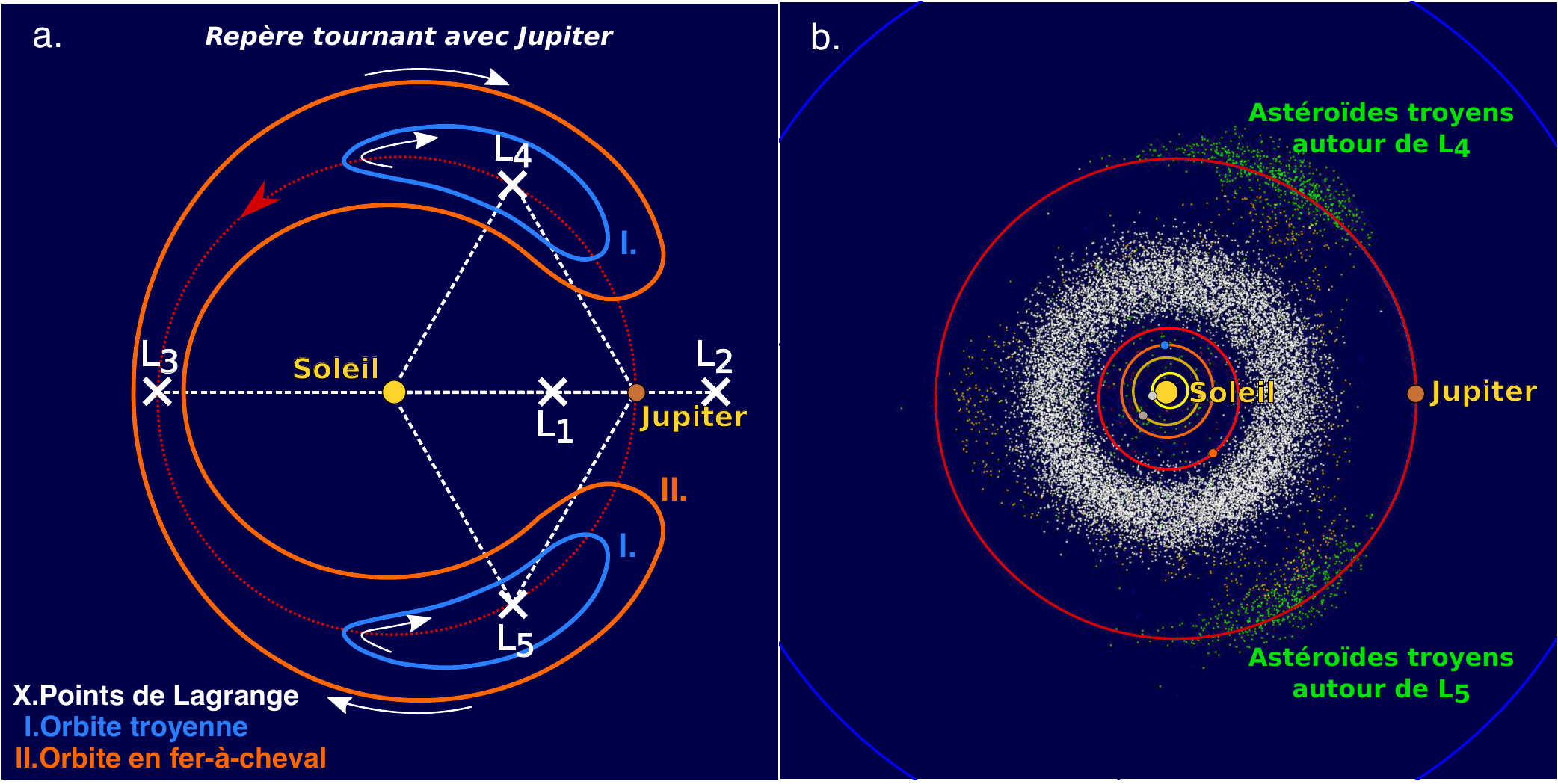}\\[0.25cm]
		\caption{\small{
			(a.) Trajectoires particuli\`eres du probl\`eme des 3 corps dans le cas d'un ast\'ero\"ide (corps de masse n\'egligeable) uniquement soumis aux forces d'attractions du Soleil et d'une plan\`ete (Jupiter) suppos\'ee en orbite circulaire. 
			Celles-ci sont repr\'esent\'ees dans un rep\`ere tournant avec Jupiter.
			Les trajectoires de type ``troyenne" (courbes bleues) passent alternativement d'une orbite interne \`a externe respectivement \`a l'orbite de Jupiter en entourant les points de Lagrange $L_4$ ou $L_5$.   
			La trajectoire de type ``fer-\`a-cheval" (courbe orange) passent \'egalement d'une orbite interne \`a externe respectivement \`a l'orbite de Jupiter mais en entourant les points de Lagrange $L_4$, $L_3$ et $L_5$. Autrement dit, elle passe contin\^ument de configurations proches du triangle \'equilat\'eral \`a la configuration align\'ee de part et d'autre du Soleil tout en ayant des rencontres proches avec la plan\`ete.
			Cette trajectoire est du m\^eme type que celles de Janus et \'Epim\'eth\'ee mais dans le cas o\`u la masse de cette derni\`ere serait n\'egligeable par rapport \`a celle de Janus (or nous savons que $m_{J}/m_{E} \simeq 3$ et donc que les masses de Janus et \'Epim\'eth\'ee sont comparables).
			(b.) Illustration du syst\`eme solaire jusqu'\`a l'orbite de Jupiter mettant en \'evidence les ast\'ero\"ides troyens  situ\'es au voisinage des position $L_4$ et $L_5$ du syst\`eme Soleil-Jupiter.
	}}
	\label{fig:QCOrbit}
	\end{center}
\end{figure}

Ces \textit{configurations align\'ees} et \textit{\'equilat\'erales} co\"incident avec les c\'el\`ebres \textit{points de Lagrange}
	\footnote{voir l'article d'Images des Math\'ematiques ``\textit{les fameux points de Lagrange}"\\ \url{http://images.math.cnrs.fr/Les-fameux-points-de-Lagrange.html}}
(not\'es respectivement $L_1$, $L_2$, $L_3$ et $L_4$, $L_5$) connus dans le cas o\`u l'un des trois corps a une  masse n\'egligeable par rapport aux deux autres.
Dans le cas particulier d'un ast\'ero\"ide uniquement soumis aux forces d'attraction du Soleil et d'une plan\`ete suppos\'ee en orbite circulaire, ce sont des trajectoires singuli\`eres car repr\'esentant des points fixes dans le r\'ef\'erentiel centr\'e sur le Soleil et tournant avec la plan\`ete (voir la figure \ref{fig:QCOrbit}.a).

\end{sloppypar}

\begin{sloppypar}

En \'enon\c cant que:
\textit{``trois masses \'etant plac\'ees non plus rigoureusement, mais \`a tr\`es-peu pr\`es dans les conditions  \'enonc\'ees [pr\'ec\'edemment], on demande si l'action r\'eciproque des masses maintiendra le syst\`eme dans cet  \'etat particulier de mouvement ou si elle tendra au contraire \`a l'en \'ecarter de plus en plus"}
	\footnote{L. Euler, ``\textit{Opera Omnia}". V.25, S.2, p. 246--257.},
Euler semblait penser que les configurations align\'ees \'etaient stables.	

Liouville d\'emontra en 1842 que les configurations align\'ees d'Euler \'etaient instables quelle que soit la valeur des masses des corps.
Par cons\'equent si la Lune avait occup\'e la position \'enonc\'ee par Euler, elle n'aurait pu s'y maintenir que pendant un temps tr\`es court.
L'ann\'ee suivante, Gascheau pose la premi\`ere pierre d'une \'eventuelle preuve de stabilit\'e de la configuration \'equilat\'erale de Lagrange sous une condition qui implique que deux masses soient suffisamment faibles par rapport \`a la troisi\`eme.
Bien plus tard, un r\'esultat de Liapounov permit d'en d\'eduire l'existence  de familles de solutions p\'eriodiques au voisinage des points d'\'equilibre $L_4$ et $L_5$.
Par ailleurs, ces solutions p\'eriodiques ont fourni une explication \`a la pr\'esence des ast\'ero\"ides ``troyens" observ\'es depuis 1906 au voisinage des positions $L_4$ et $L_5$ du syst\`eme Soleil-Jupiter (voir la figure \ref{fig:QCOrbit}.a--b).

Contrairement aux configurations d'Euler et de Lagrange, ces orbites p\'eriodiques \textit{troyennes}  ne d\'ecrivent pas exactement des ellipses autour du Soleil mais en sont tr\`es proches: ce sont de petites d\'eformations p\'eriodiques de la configuration \'equilat\'erale.
Ainsi, ce sont des solutions perp\'etuellement stables du probl\`eme des 3 corps, en r\'esonance co-orbitale et qui poss\`edent notamment la propri\'et\'e d'\'echange d'orbites.
En effet, comme illustr\'e sur la figure \ref{fig:QCOrbit}.a, l'ast\'ero\"ide troyen passe alternativement d'orbite externe \`a interne par rapport \`a celle de Jupiter (qui reste la m\^eme).
Par contre, le mouvement relatif de l'ast\'ero\"ide et de la plan\`ete entoure uniquement l'un des points de Lagrange $L_4$ ou $L_5$: il ne d\'ecrit donc pas une trajectoire en fer-\`a-cheval comme celle de Janus et \'Epim\'eth\'ee autour de Saturne.

\end{sloppypar}

\begin{sloppypar}

\section{Notre probl\`eme : la ``stabilit\'e des fers-\`a-cheval"}

Revenons maintenant au cas des lunes co-orbitantes de Saturne.
Nous souhaitions d\'emontrer l'existence de solutions perp\'etuellement stables associ\'ees aux trajectoires en fer-\`a-cheval  dans le cadre du probl\`eme des 3 corps
\footnote{Plus pr\'ecis\'ement dans le cadre du probl\`eme des 3 corps o\`u deux masses comparables gravitent dans un m\^eme plan autour d'un corps central beaucoup plus massif.}.
Il s'agissait donc d'\'etudier une version simplifi\'ee du probl\`eme r\'eel en n\'egligeant les autres forces en action:  l'influence gravitationnelle d'autres satellites de Saturne, l'effet des anneaux, l'aplatissement de Saturne, les effets de mar\'ees, etc...

Nous recherchions des trajectoires en r\'esonance co-orbitale qui, comme dans la figure \ref{fig:QCOrbit}.a dans le cas d'un syst\`eme Soleil-plan\`ete-ast\'ero\"ide (courbe orange),   ``entourent" les trois configurations d'\'equilibre  $L_4$, $L_3$ et $L_5$, c'est-\`a-dire qui passent continument de configurations proches du triangle \'equilat\'eral \`a la configuration align\'ee de part et d'autre de Saturne tout en ayant des rencontres proches entre les deux lunes.
Or, contrairement aux orbites p\'eriodiques troyennes d\'ecrites pr\'ec\'edemment, ce type de configuration ne correspond pas \`a de petites d\'eformations d'une configuration d'\'equilibre.
Dans ce cas, les math\'ematiques requises ne sont plus les m\^emes.
Ainsi, bien qu'elle ait \'et\'e pressentie au d\'ebut du vingti\`eme si\`ecle  par l'astronome E. Brown\footnote{E. W. Brown ``\textit{On a New Family of Periodic Orbits in the Problem of Three Bodies}", 1911, Monthly Notices of the Royal Astronomical Society, vol.71, Issue 5, p.438--454.} dans une \'etude tr\`es simplifi\'ee du probl\`eme des 3 corps, l'existence de trajectoires en fer-\`a-cheval perp\'etuellement stables n'avait jamais \'et\'e \'etablie avant notre r\'esultat.
Comme nous le verrons dans les paragraphes suivant, ce sont des math\'ematiques des ann\'ees 90 qui nous ont permis de le prouver.

\end{sloppypar}

\section{Une autre mani\`ere de trouver des solutions: la th\'eorie KAM}

\begin{sloppypar}
Au d\'ebut du vingti\`eme si\`ecle, Poincar\'e fournit un r\'esultat crucial qui confirma la complexit\'e des solutions du probl\`eme des 3 corps en d\'emontrant sa \textit{non-int\'egrabilit\'e}.
En un certain sens, cela revient \`a montrer que les solutions du probl\`eme des 3 corps sont qualitativement diff\'erentes des mouvements kepleriens qui sont \textit{quasi-periodiques} (c'est-\`a-dire comme la superposition des mouvements p\'eriodiques des diff\'erents corps c\'elestes).
En effet, Poincar\'e exhiba une multiplicit\'e de solutions tr\`es complexes, d\^ites \textit{chaotiques}
	\footnote{Pour aller plus loin, voir la partie d\'edi\'ee \`a la notion math\'ematique de chaos (\url{http://ciel.mmi-lyon.fr/dynamique-chaotique/}) de l'exposition virtuelle ``\textit{les math\'ematiques du ciel}".}.
N\'eanmoins il n'excluait pas la co-existence de trajectoires simples (p\'eriodiques ou quasi-p\'eriodiques) au milieu de ces trajectoires compliqu\'ees.

Une avanc\'ee spectaculaire fut accomplie au cours ann\'ees 50-60 avec l'av\`enement de la th\'eorie KAM du nom des trois math\'ematiciens Kolmogorov, Arnol'd et Moser.
Dans un contexte plus g\'en\'eral que le probl\`eme des 3 corps, cette th\'eorie a pour objet l'\'etude de la persistance de mouvements quasi-p\'eriodiques sous l'effet d'une perturbation en combinant syst\`emes dynamiques, arithm\'etique et probabilit\'es.
Ainsi, si l'on conna\^it une \textit{approximation int\'egrable}, c'est-\`a-dire une approximation du probl\`eme consid\'er\'e par un autre mais dont les solutions sont compos\'ees uniquement de mouvements quasi-p\'eriodiques, la th\'eorie KAM garantit qu'en dehors de certains cas complexes, l'effet d'une perturbation modifie la forme des solutions de l'approximation int\'egrable sans toutefois d\'etruire leur caract\`ere quasi-p\'eriodique
	\footnote{Voir l'article r\'ecent de B. Gr\'ebert dans Images des Math\'ematiques \url{http://images.math.cnrs.fr/Du-ressort-a-l-univers-quand-les-resonances-se-jouent-des-pronostics}}.
Plus pr\'ecis\'ement, cette th\'eorie s'applique pour des trajectoires fortement non-r\'esonantes (l'arithm\'etique intervient alors pour quantifier la notion de ``distance" aux r\'esonances) qui forment un ``gros" ensemble parmi l'ensemble des solutions du probl\`eme (c'est alors la notion de ``mesure" de cet ensemble qui fait intervenir la th\'eorie des probabilit\'es).
Il s'agissait d'un r\'esultat de stabilit\'e retentissant tant l'ampleur des difficult\'es math\'ematiques \`a surmonter \'etait importante.

En 1963, Arnol'd
	\footnote{V.I. Arnol'd ``\textit{Small Denominators and Problems of Stability of Motion in Classical and Celestial Mechanics}", 1963, Russian Mathematical Surveys, vol.18, p.85--191.}
appliqua ce r\'esultat au probl\`eme des 3 corps pour deux corps de masses comparables gravitant dans un m\^eme plan autour d'un troisi\`eme beaucoup plus massif, et d\'emontra le th\'eor\`eme suivant:
\textbf{\textit{``si les masses de Jupiter et Saturne avaient \'et\'e suffisamment petites par rapport \`a celle du Soleil, pour beaucoup (dans un sens math\'ematique pr\'ecis) de conditions initiales, le mouvement de ces plan\`etes aurait \'et\'e quasi-p\'eriodique"}}.
Plus r\'ecemment, ce r\'esultat a \'et\'e \'etendu  par M. Herman et J. F\'ejoz au cas de $N$ corps (plan\`etes) de masses comparables gravitant dans l'espace autour d'un corps central (le Soleil) beaucoup plus massif.

Bien qu'il ne soit pas adapt\'e au cas r\'eel puisqu'il n\'ecessite des masses incomparablement plus petites que les masses des plan\`etes du syst\`eme solaire, ce th\'eor\`eme garantit l'existence de solutions quasi-p\'eriodiques pour deux (ou $N$) plan\`etes en mouvement autour du Soleil sans que l'effet de leurs interactions mutuelles puisse d\'estabiliser leurs trajectoires en provoquant des \'ejections, des collisions ou des r\'e-arrangements entre les orbites. 
Nous pouvons cependant pr\'eciser que des \'etudes num\'eriques men\'ees dans le cas des masses du Soleil, de Jupiter et de Saturne fournissent des solutions extr\^emement stables sur des temps comparables \`a l'\^age du syst\`eme solaire et qui semblent quasi-p\'eriodiques, ce qui s'accorde avec les pr\'edictions de la th\'eorie KAM.
Malheureusement \`a ce jour, aucune preuve math\'ematique de ce comportement n'a \'et\'e d\'emontr\'ee.

\end{sloppypar}

\section{Th\'eorie KAM et ``fers-\`a-cheval"}

\begin{figure}[H]
	\begin{center}
		\includegraphics[scale=0.575]{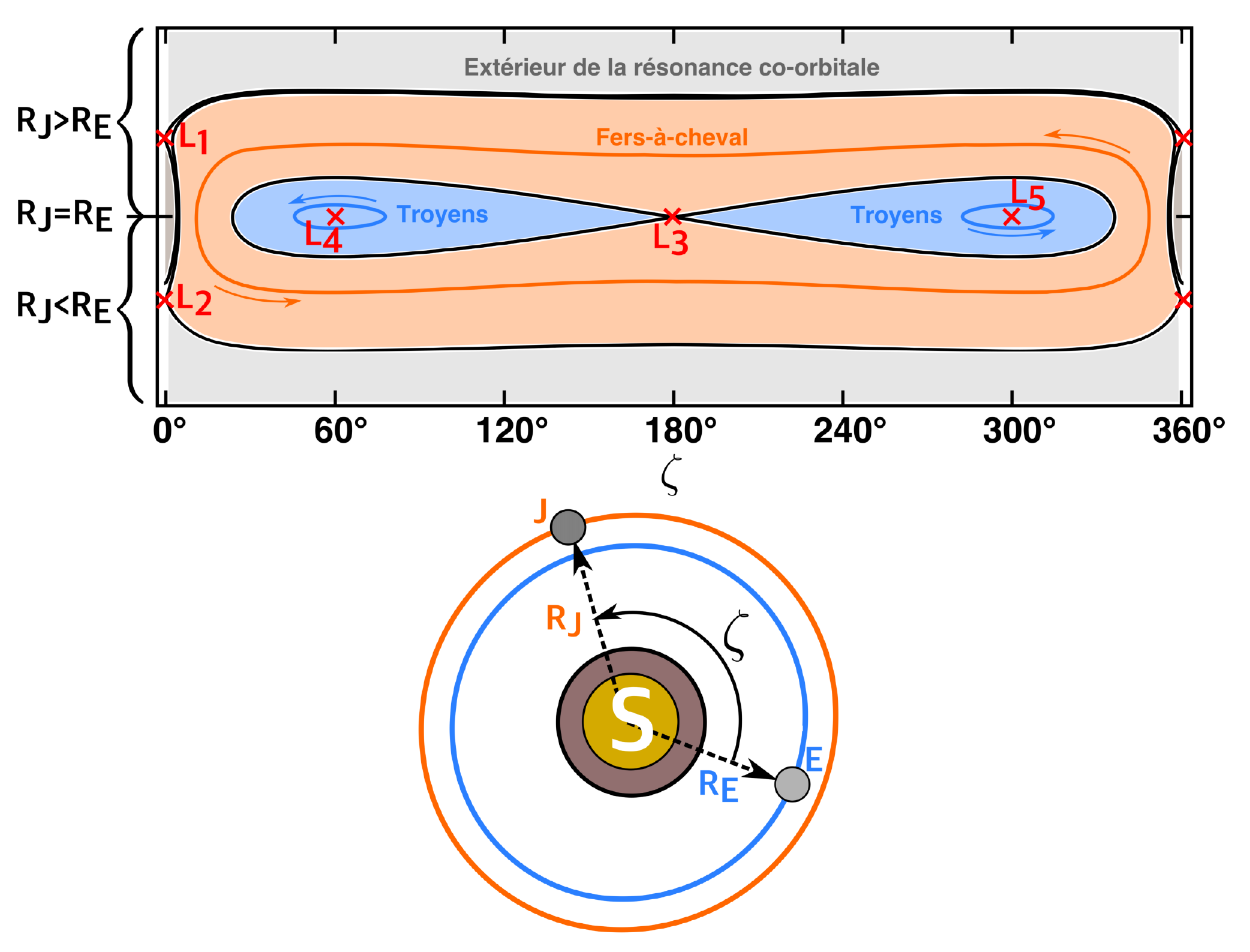}\\[0.25cm]
		\caption{\small{
		Dynamiques fournies par notre approximation int\'egrable dans les variables r\'esonantes: la diff\'erence entre les rayons des orbites circulaires des deux lunes, $R_E-R_J$ et l'angle $\zeta = \widehat{\bf{JSE}}$. Voir l'encadr\'e ``description du mod\`ele" pour plus de d\'etails.
		}}
\label{fig:Phase}
	\end{center}
\end{figure}


\begin{sloppypar}

Poursuivant l'id\'ee d'Arnol'd, c'est la th\'eorie KAM que nous appliquons au cas du syst\`eme Saturne-Janus-\'Epim\'eth\'ee mais dans un contexte  plus d\'elicat.
En effet, alors que la th\'eorie d'Arnol'd concerne des mouvements kepleriens fortement non-r\'esonants (voir plus haut), nous consid\'erons une situation r\'esonante puisque les deux lunes sont co-orbitantes. 
En outre, alors qu'Arnol'd pouvait s'appuyer sur les mouvements kepleriens comme mouvements de r\'ef\'erence, notre cas n\'ecessitait de construire une approximation int\'egrable dont les solutions d\'ecrivent enti\`erement la topologie de la r\'esonance co-orbitale: les configurations  d'\'equilibre $L_3$, $L_4$ et $L_5$, la dynamique des orbites troyennes qui oscillent autour de $L_4$ et $L_5$ et la dynamique des trajectoires en fer-\`a-cheval qui entourent $L_3$, $L_4$ et $L_5$.
La figure 7 illustre les diff\'erents types de trajectoires associ\'es \`a ce mod\`ele int\'egrable. 

\end{sloppypar}

 \medskip

\noindent\fbox{\parbox{\linewidth\fboxrule\fboxsep}{\textbf{Description du mod\`ele.}\\
Le long d'une solution de notre mod\`ele, Janus et \'Epim\'eth\'ee ont des trajectoires  circulaires de rayons $R_J$  et $R_E$ autour de Saturne.
Dans la figure \ref{fig:Phase}, l'axe des ordonn\'ees repr\'esente l'angle $\zeta$ form\'e par Janus, Saturne et \'Epim\'eth\'ee tandis que l'axe  des abscisses est associ\'e aux d\'eformations de leur rayon. 

Lorsque les lunes sont align\'ees telle que $\zeta=0$, nous retrouvons les deux configurations align\'ees $L_1$ et $L_2$. 
De celles-ci apparaissent des ``s\'eparatrices" (les courbes noires passant par $L_1$ et $L_2$) qui divisent le domaine en deux parties o\`u les trajectoires ont des comportements bien distincts: 
\begin{itemize}
\item l'ext\'erieur non-r\'esonant (r\'egions grises) o\`u la diff\'erence entre les rayons (et par cons\'equent entre les p\'eriodes orbitales) est importante, 
\item l'int\'erieur correspondant \`a la ``zone" de r\'esonance co-orbitale  (mentionn\'ee dans l'encadr\'e d\'etaillant la notion de r\'esonance) o\`u les \'echanges d'orbites se produisent (les trajectoires passent toutes par $R_J = R_E$).
\end{itemize}   
   
Lorsque les rayons sont \'egaux et que $\zeta$ est \'egal $60$, $300$ et $180$ degr\'es, nous retrouvons respectivement les configurations \'equilat\'erales $L_4$, $L_5$ et la configuration align\'ee $L_3$.

Enfin, une s\'eparatrice passant par $L_3$ divise la ``zone" de r\'esonance en deux  dynamiques diff\'erentes:
\begin{itemize}
\item les trajectoires ``troyennes" (les domaines bleus) le long desquelles les lunes \'echangent leurs orbites en maintenant une configuration proche de celle \'equilat\'erale  (l'angle $\zeta$ oscille autour de $60$ ou de $300$  degr\'es comme pour les deux courbes bleues),
\item les trajectoires en fer-\`a-cheval (domaine rouge) le long desquelles, l'angle $\zeta$ oscille autour de $180$ degr\'es avec une tr\`es grande amplitude jusqu'\`a ce que les lunes se rapprochent et \'echangent leurs orbites (lorsque $|\zeta|$ est petit; voir la courbe rouge).
\end{itemize}
}}

\medskip

\begin{sloppypar}
D'autre part, n'ayant pas d'expression explicite des solutions de l'approximation int\'egrable, il n'a pas \'et\'e possible d'utiliser les versions classiques de la th\'eorie KAM.
De ce fait, nous avons  \'et\'e contraints d'utiliser une version plus r\'ecente de cette th\'eorie, d\'evelopp\'ee
\footnote{J. P\"oschel ``\textit{A KAM-theorem for some nonlinear partial differential equations}", 1996, Annali della Scuola Normale Superiore di Pisa, Classe di
Scienze, vol.23, p.119--148.}
 par J. P\"oschel dans les ann\'ees 90, o\`u les conditions requises sont particuli\`erement faibles.

Une fois ces difficult\'es surmont\'ees, nous avons \'etabli le th\'eor\`eme
	\footnote{L. Niederman, A. Pousse, P. Robutel ``\textit{On the co-orbital motion in the three-body problem: existence of quasi-periodic horseshoe-shaped orbits}", 2018. Preprint disponible sur Arxiv: \url{https://arxiv.org/abs/1806.07262}}
suivant:
\textbf{\textit{``en supposant que les masses de Janus et \'Epim\'eth\'ee soient suffisamment petites par rapport \`a celle de Saturne, il existe des conditions initiales telles que leur ballet gravitationnel continuera ind\'efiniment.}}
\end{sloppypar}

\end{document}